\documentclass{article}
\usepackage{graphicx}
\usepackage{epstopdf}
\usepackage{amsmath}
\usepackage{amssymb}
\usepackage{bm}
\usepackage{bbm}
\usepackage{enumerate}
\usepackage[T1]{fontenc}
\usepackage[latin1]{inputenc}
\usepackage{hyperref}
\usepackage{color}
\usepackage{curves}
\usepackage{tikz}
\usepackage[hang]{caption}
\usepackage{ntheorem}
\usepackage{float}
\usepackage{mathrsfs}

\title{Water transport on infinite graphs}
\author{
\makebox[\textwidth][c]{
	\begin{minipage}[t]{\textwidth}
		\begingroup
		\begin{minipage}[t]{0.5\textwidth}
			Olle Häggström \thanks{Research supported by grants
					from the Swedish Research Council and from the Knut and Alice
					Wallenberg Foundation}\\\normalsize Chalmers University of Technology
		\end{minipage}
		\hfill
		\begin{minipage}[t]{0.5\textwidth}
			Timo Hirscher
				\thanks{Research supported by grants from the Swedish Research Council and the Royal Swedish Academy of Sciences}\\\normalsize Goethe-University Frankfurt am Main
		\end{minipage}
		\endgroup
\end{minipage}}}

\theoremstyle{break}
\newtheorem{theorem}{Theorem}[section]
\newtheorem{corollary}{Corollary}[section]
\newtheorem{lemma}{Lemma}[section]

\theorembodyfont{\upshape}
\newtheorem{definition}{Definition}
\newtheorem*{remark}{Remark}
\newtheorem*{example}{Examples}

\makeatletter
\let\c@proposition\c@theorem
\let\c@lemma\c@theorem
\let\c@corollary\c@theorem
\makeatother

\newenvironment{proof}{\noindent{\sc Proof:}}{\vspace{-0.5cm}~\hfill $\square$\vspace{0.5cm}}
\newenvironment{nproof}[1]{\noindent{\sc Proof #1:}}{\vspace{-1em}~\hfill $\square$\vspace{2em}}

\newcommand\N{\mathbb{N}}
\newcommand\R{\mathbb{R}}
\newcommand\Z{\mathbb{Z}}
\newcommand\E{\mathbb{E}\,}
\newcommand\Prob{\mathbb{P}}
\renewcommand\epsilon{\varepsilon}
\renewcommand\phi{\varphi}
\DeclareMathOperator\supp{supp}

\definecolor{darkblue}{rgb}{0,0,.5}
\hypersetup{colorlinks=true, breaklinks=true, linkcolor=blue, 
citecolor=darkblue, menucolor=blue, urlcolor=blue}

\begin{document}
\newpage
\maketitle
\begin{abstract}
    If the nodes of a graph are considered to be identical barrels -- featuring different water
    levels -- and the edges to be (locked) water-filled pipes in between the barrels,
    consider the optimization problem of how much the water level in a fixed barrel can be raised
    with no pumps available, i.e.\ by opening and closing the locks in an elaborate succession.
    This problem originated from the analysis of an opinion formation process, the so-called Deffuant model. We consider i.i.d.\ random initial water levels and ask whether the supremum of achievable levels at a given node has a degenerate distribution. This turns out to be the case for all infinite connected quasi-transitive graphs with exactly one exception: the two-sided infinite path.
\end{abstract}


\section{Introduction}
\subsubsection*{Motivation and model}

\begin{figure}[ht]
     \centering
     \includegraphics[scale=0.9]{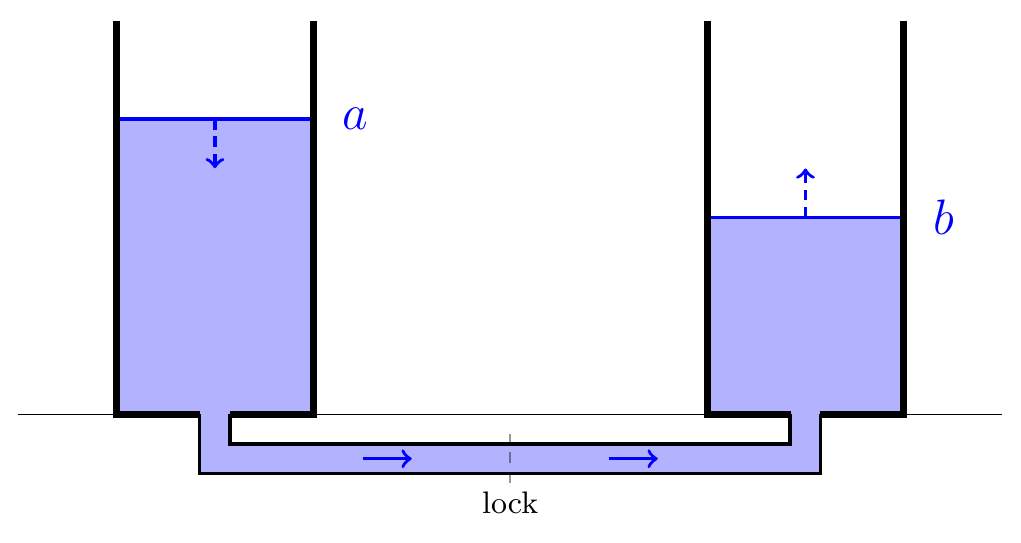}
     \caption{Leveling water stages after just having opened a lock. \label{barrels}}
\end{figure}

Imagine a plane on which rainwater is collected in identical rain barrels, some of which are connected through
pipes (that are already water-filled). All the pipes feature locks that are normally closed. If a lock is opened,
the contents of the two barrels which are connected via this pipe start to level, see Figure \ref{barrels}.
If one waits long enough, the water levels in the two barrels will be exactly the same, namely lie at the average
$\tfrac{a+b}{2}$ of the two water levels ($a$ and $b$) before the pipe was unlocked.

After a rainy night in which the barrels accumulated various amounts of precipitation we might be interested
in maximizing the water level in one fixed barrel by opening and closing some of the locks in carefully chosen
order.

\vspace*{1em}
\noindent 
In order to mathematically model the setting, consider an undirected graph $G=(V,E)$, which is either finite or infinite with bounded maximum degree. Furthermore, we can assume without loss of generality
that $G$ is connected and simple, that means having neither loops nor multiple edges. Every vertex is
understood to represent one of the barrels and the pipes correspond to the edges in the graph. The
barrels themselves are considered to be identical, having a fixed capacity, say $1$ for simplicity.

\subsubsection*{Making a move}
Given some initial profile $\{\eta_0(u)\}_{u\in V}\in[0,1]^V$, the system is considered to evolve in discrete
time and in each round we can open one of the locked pipes and
transport water from the fuller barrel into the emptier one. If we stop early, the two levels might not have
completely balanced out giving rise to the following {\em update rule} for the water profile: If in round $k$ the pipe
$e_k=\langle x,y\rangle$ connecting the two barrels at sites $x$ and $y$, with levels $\eta_{k-1}(x)=a$ and
$\eta_{k-1}(y)=b$ respectively, is opened and closed after a certain period of time, we get

 \begin{equation}\label{update}
         \begin{array}{rcl}\eta_k(x) &\!=\!& a+\mu_k\,(b-a)\\
                           \eta_k(y) &\!=\!& b+\mu_k\,(a-b)
         \end{array}
 \end{equation}
for some $\mu_k\in(0,\tfrac12]$, which we assume can be chosen freely by an appropriate choice of how long the pipe
is left open. All other levels stay unchanged, i.e.\ $\eta_k(w)=\eta_{k-1}(w)$ for all $w\in V\setminus\{x,y\}$.

Such a move can be described by the tuple $(e_k,\mu_k)\in E\times(0,\frac12]$, consisting of chosen edge and transferred fraction. A {\em finite sequence of moves} is accordingly the concatenation of
$(e_1,\mu_1),\dots,(e_T,\mu_T)$, in chronological order, and leads to a final water profile $\{\eta_T(u)\}_{u\in V}$.
The optimization problem is to maximize $\eta_T(v)$ over $T\in\N_0$ and move sequences. The quantity of interest is defined as follows:

\begin{definition}\label{kappa}
	For a graph $G=(V,E)$, an initial water profile $\{\eta_0(u)\}_{u\in V}$ and a fixed {\em target vertex} $v\in V$, let $\kappa_G(v)$ be defined as the supremum over all water levels that are
	achievable at $v$ with finite move sequences, i.e.\
	$$\kappa_G(v):=\sup\{r\geq 0:\, \text{there exists }T\in\N_0 \text{ and a move sequence s.t.\ }\eta_T(v)=r\}.$$
\end{definition}

\subsubsection*{Connections to other concepts and synopsis}
Some readers, familiar with mathematical models for social interaction processes, might note that (\ref{update})
basically looks like the update rule in the opinion formation process given
by the so-called {\em Deffuant model} for consensus formation in social networks (as described in the
introduction of \cite{Deffuant}), leaving aside the fact that here $\mu$ can change from update to update and the bounded confidence restriction is omitted. This similarity, however, is no coincidence: The situation described in the context above arises naturally in the analysis of the Deffuant model when dealing with the question of how extreme the opinion of a fixed agent can get (given an initial opinion profile on a specified network graph) if the interactions take place appropriately.

In order to tackle this question, Häggström \cite{ShareDrink} invented a non-random pairwise averaging
procedure, which he proposed to call {\em Sharing a drink} (SAD). This process -- which is the main
focus of the preparatory Section \ref{sec2} -- was originally considered only on the (two-sided)
infinite path, i.e.\ the graph $G=(V,E)$ with $V=\Z$ and $E=\{\langle v, v+1\rangle,\;v\in\Z\}$, but
can immediately be generalized to any graph (see Definition \ref{SAD}) and is dual to the water
transport described above in a sense to be made precise in Lemma \ref{dual}.

In Section \ref{infinite}, we consider the initial water levels to be i.i.d.\ random variables, uniformly distributed on $[0,1]$, and detect a remarkable change of behavior: On the infinite path $\Z$, the highest achievable water level at a fixed vertex depends on the initial profile in the sense that it has a nondegenerate distribution, just like on any finite graph.
If the infinite graph contains a neighbor-rich half-line (see Definition \ref{hl}), however, this
dependence becomes degenerate: For any vertex $v\in V$, the value $\kappa_G(v)$ almost surely equals
1 (the essential supremum of the marginal distribution).

This fact makes the two-sided infinite path quite unique:
It constitutes the only exception among all infinite connected quasi-transitive graphs, to the effect
that $\kappa_G(v)$ is a nondegenerate random variable -- cf.\ Corollary \ref{main}. Furthermore, this
result is interesting also in view of the connection to opinion formation: It indicates one reason why
the analysis of the Deffuant model on $\Z^d$ with $d\geq2$ has so far encountered more resistance as compared to the case of $d=1$ (see the concluding Remark part (\ref{concl}) for details).

\section{Connection to the SAD-procedure}\label{sec2}

Let us first recall from \cite{ShareDrink} the formal definition of the SAD-procedure:
\begin{definition}\label{SAD}
For a graph $G=(V,E)$ and some fixed vertex $v\in V$, define $\{\xi_0(u)\}_{u\in V}$ by setting
$$\xi_0(u)=\delta_v(u):=\begin{cases}1\quad\text{for } u=v\\ 0 \quad\text{for } u\neq v.\end{cases}$$
In each time step, an edge $\langle x,y\rangle$ is chosen and the profile $\{\xi_0(u)\}_{u\in V}$ updated
according to the rule (\ref{update}) with $\{\xi_k(u)\}_{u\in V}$ in place of $\{\eta_k(u)\}_{u\in V}$.
One can interpret this process as a full glass of water initially placed at vertex $v$ (all other glasses
being empty), which is then repeatedly shared among neighboring vertices by (each time step) choosing a pair
of neighbors and pouring a $\mu_k$-fraction of the difference from the glass containing more water into the one
containing less. Let us refer to this interaction process as {\em Sharing a drink (SAD)}.
\end{definition}

Just as in \cite{ShareDrink}, the SAD-procedure can be used to describe the composition of the contents
in the water barrels after finitely many rounds of opening and closing pipe locks. The following lemma
corresponds to Lemma 3.1 in \cite{ShareDrink}, albeit in a more general graph context. Since the two dual processes (water transport and SAD) evolve
in discrete time in our setting, its proof simplifies somewhat.

\begin{lemma}\label{dual}
Consider an initial profile of water levels $\{\eta_0(u)\}_{u\in V}$ on a graph $G=(V,E)$ and fix a
vertex $v\in V$. For $T\in\N_0$ define the SAD-procedure that starts with $\xi_0=\delta_v$
(see Definition \ref{SAD}) and is dual to the chosen move sequence in the water transport problem in the 
following sense: If in round $k\in\{1,\dots,T\}$ the water profile is updated according to (\ref{update}), the update in the SAD-profile at
time $T-k\in\{0,\dots,T-1\}$ takes place along the same edge and with the same choice of $\mu_k$. Then we get
\begin{equation}\label{convcomb}\eta_T(v)=\sum_{u\in V}\xi_T(u)\,\eta_0(u).\end{equation}
\end{lemma}

\begin{proof}
We prove the statement by induction on $T$. For $T=0$, the statement is trivial and there is nothing to show.
For the induction step fix $T\in\N$ and assume the first pipe opened to be $e=\langle x,y\rangle$. According to 
(\ref{update}) we get
$$\eta_1(u)=\begin{cases}\eta_0(u)& \text{if }u\notin\{x,y\}\\
                         (1-\mu_1)\,\eta_0(x)+\mu_1\,\eta_0(y)&\text{if }u=x\\
                         (1-\mu_1)\,\eta_0(y)+\mu_1\,\eta_0(x)&\text{if }u=y.
                         \end{cases}$$
Let us consider $\{\eta_1(u)\}_{u\in V}$ as some initial profile $\{\eta'_0(u)\}_{u\in V}$. By induction
hypothesis we get
\begin{align*}
\eta'_{T-1}(v)&=\sum_{u\in V}\xi'_{T-1}(u)\,\eta'_0(u)\\
            &=\sum_{u\in V\setminus\{x,y\}}\xi'_{T-1}(u)\,\eta_0(u)
            +\Big((1-\mu_1)\,\xi'_{T-1}(x)+\mu_1\,\xi'_{T-1}(y)\Big)\,\eta_0(x)\\
            &\quad+\Big((1-\mu_1)\,\xi'_{T-1}(y)+\mu_1\,\xi'_{T-1}(x)\Big)\,\eta_0(y),
\end{align*}
where $\eta'_{T-1}(v)=\eta_T(v)$ and $\{\xi'_k(u)\}_{u\in V}$, $0\leq k\leq T-1$, are the intermediate water
profiles of the SAD-procedure corresponding to the move sequence after round $1$. As by definition the original
SAD-procedure arises from the shortened one by adding an update at time $T-1$ along edge $e$ with parameter $\mu_1$,
we find $\xi_k(u)=\xi'_k(u)$ for all $k\in\{0,\dots,T-1\}$ and $u\in V$ as well as
$$\xi_T(u)=\begin{cases}\xi_{T-1}(u)=\xi'_{T-1}(u)& \text{if }u\notin\{x,y\}\\
                         (1-\mu_1)\,\xi_{T-1}(x)+\mu_1\,\xi_{T-1}(y)&\text{if }u=x\\
                         (1-\mu_1)\,\xi_{T-1}(y)+\mu_1\,\xi_{T-1}(x)&\text{if }u=y,
                         \end{cases}$$
which establishes the claim.
\end{proof}

\vspace*{1em}
Note that, although we are going to consider random initial profiles of
water levels, the statement
of Lemma \ref{dual} deals with a deterministic duality that does not involve any randomness (once the initial
profile and the move sequence are fixed).

Let us prepare one more auxiliary result, which turns out to follow readily from the energy argument used in the proof of Thm.\ 2.3 in \cite{ShareDrink}:
\begin{lemma}\label{evenout}
Given an initial profile of water levels $\{\eta_0(u)\}_{u\in V}$ on a graph $G=(V,E)$, fix a finite set $A\subseteq V$
and a set $E_A\subseteq E$ of edges inside $A$ that connects $A$.
If we open the pipes in $E_A$ -- and no others -- in repetitive sweeps for times long enough
such that $\mu_k\geq\epsilon$ for some fixed $\epsilon>0$ in each round (cf.\ (\ref{update})),
then the water levels at vertices in $A$ approach a balanced average, i.e.\ 
converge to the value $\tfrac{1}{|A|} \sum_{u\in A}\eta_0(u)$. For all $v\in A$, the corresponding dual SAD-profiles started with
$\xi_0=\delta_v$ converge uniformly to $\tfrac{1}{|A|}\,\mathbbm{1}_A$.
\end{lemma}

\begin{proof}
 Let us define the energy after round $k$ inside $A$ by
 $$W_k(A)=\sum_{u\in A}\big(\eta_k(u)\big)^2.$$
 A short calculation reveals that an update of the form (\ref{update}) reduces the energy by
 $2\mu_k^{\,2}\,(b-a)^2$, where the updated water levels were $a$ and $b$ respectively. If $\mu_k$ is bounded away from $0$,
 the fact that $W_k(A)\geq0$ for all $k$ entails that the difference in water levels $|b-a|$ before a pipe is opened
 can be larger than any fixed positive value only finitely many times. In effect, since any pipe in $E_A$ is opened
 repetitively we must have $|\eta_k(u)-\eta_k(v)|\to 0$ as $k\to\infty$ for all edges $\langle u,v\rangle\in E_A$.
 As the updates are average preserving, the first part of the claim follows from the fact that $E_A$ connects $A$.
 
 The second part of the lemma follows by applying the same argument to the dual SAD-procedure.
 \end{proof}
 \vspace*{1em}
 
 To round off these preliminary considerations, let us collect some properties of SAD-profiles derived in $\cite{ShareDrink}$ into a single lemma for convenience.
 
 \begin{lemma}\label{collection}
  Consider the SAD-procedure on a path, started in vertex $v$, i.e.\ with $\xi_0(u)=\delta_v(u),\ u\in V$.
  \begin{enumerate}[(a)]
   \item The SAD-profiles achievable on paths are all unimodal.
   \item If the vertex $v$ only shares the water to one side, it will remain a mode of the SAD-profile.
   \item The supremum over all achievable SAD-profiles started with $\delta_v$ at another vertex $w$ equals $\tfrac{1}{d+1}$,
         where $d$ is the graph distance between $v$ and $w$.
 \end{enumerate}
 \end{lemma}	
 
  The results in \cite{ShareDrink} actually all deal with the two-sided infinite path, but it is evident how the arguments
  used immediately transfer to finite paths. Part (a) hereby corresponds to Lemma 2.2 in \cite{ShareDrink}, part
  (b) to Lemma 2.1 and part (c) to Thm.\ 2.3.
  The argument Häggström \cite{ShareDrink} used to prove the statement in (c) for the two-sided infinite path can in fact
  be generalized to prove the result for trees without much effort, as was done by Shang (see Prop.\ 6 in \cite{Shang}).

\section{In terms of water transport, \texorpdfstring{$\Z$}{Z} behaves like a finite graph}\label{infinite}

In this section, we want to analyze the water transport problem given random
initial water levels. More precisely, we will consider connected, simple graphs $G=(V,E)$ with bounded degree and i.i.d.\ $\mathrm{unif}(0,1)$ initial water levels $\{\eta_0(u)\}_{u\in V}$. The supremum $\kappa_G(v)$ of achievable water levels at a fixed target vertex $v\in V$ (cf.\ Def.\ \ref{kappa}) naturally depends on the initial water levels, which makes it a random variable as well.
The question we want to address is in which cases the distribution of $\kappa(G)$ is degenerate.

Let us back away from $\kappa_G(v)$ for a moment and briefly discuss the global average of water
levels across the graph (in our i.i.d.\ $\mathrm{unif}(0,1)$ setting), which, as a side note, will not
change with time as updates are average-preserving, cf.\ \eqref{update}.
Here the picture is clear-cut:
On any finite graph, the average is a nondegenerate random variable, while (defined as the limit
of averages along a given sequence of nested subsets of the vertex set, which finally include every
fixed vertex) on any infinite graph it almost surely equals the expected value $\frac12$, according to
the strong law of large numbers.

For $\kappa_G(v)$, in contrast, the regime of nondegeneracy extends a bit into the realm
of infinite graphs, but not much. The most clear-cut statement of this -- Theorem \ref{qtgraphs} below -- is for the class of quasi-transitive graphs.

\begin{definition}\label{quasitrans}
	Let $G=(V,E)$ be a simple graph. A bijection $f:V\to V$ with the property that $\langle f(u),f(v)\rangle\in E$
	if and only if $\langle u,v\rangle\in E$ is called a {\em graph automorphism}. $G$ is said to be {\em (vertex-)
		transitive} if for any two vertices $u,v\in V$ there exists a graph automorphism $f$ that maps $u$ on $v$, i.e.\
	$f(u)=v$.
	If the vertex set $V$ can be partitioned into finitely many classes such that for any two vertices $u,v$
	belonging to the same class there exists a graph automorphism that maps $u$ on $v$, the graph $G$ is called
	{\em quasi-transitive}.
\end{definition}

Note that the idea of quasi-transitivity becomes nontrivial only for infinite graphs, as all finite
graphs are quasi-transitive by definition. In what follows, we will denote the support of the distribution $\mathcal{L}(X)$ of a random variable $X$ by
$$\mathrm{supp}(X):=\big\{r\in\R,\;\forall\,\epsilon>0:\ \Prob\big(|X-r|<\epsilon\big)>0\big\}.$$

\begin{theorem}\label{qtgraphs}
	Consider an infinite, connected, quasi-transitive graph $G=(V,E)$ and the initial water levels to be i.i.d.\
	$\mathrm{unif}(0,1)$. Let $v\in V$ be a fixed vertex of the graph.
	If $G$ is the two-sided infinite path $\Z$, i.e.\ $V=\Z$,  $E=\{\langle u, u+1\rangle,\;u\in\Z\}$, then $$\Prob\big(\tfrac12<\kappa_\Z(v)<1\big)=1\quad \text{and}\quad \{\tfrac12,1\}\subseteq\mathrm{supp}\big(\kappa_\Z(v)\big).$$
	If $G$ is not the two-sided infinite path, then $\kappa_G(v)=1$ almost surely.
\end{theorem}

With the intention of phrasing the special role of $\Z$ in the class of infinite (connected) quasi-transitive graphs more concisely, we can formulate the following corollary of the above theorem:

\begin{corollary}\label{main}
	Given i.i.d.\ $\mathrm{unif}(0,1)$ initial water levels, the only infinite, connected, quasi-transitive graph $G$, for which $\kappa_G(v)$ is
	not deterministic is the two-sided infinite path $\Z$.
\end{corollary}

Before setting out for proving Theorem \ref{qtgraphs}, let us get acquainted with the optimization problem of rising the water level at a fixed vertex given a random initial profile of water levels by looking at two toy examples.

\begin{example}
\begin{enumerate}[(a)]
	\item 
Consider the simplest nontrivial graph $G$, consisting of a single edge, and let the initial water levels be given by two random
variables $U_1$ and $U_2$. Trivially, we have
$$U_1\leq \kappa_G(1)\leq \max\{U_1, U_2\}.$$

If we assume $U_1$ and $U_2$ to be independent and uniformly distributed on $[0,1]$, a short calculation
reveals the distribution function
$$F_{\kappa_G(1)}(x)=\begin{cases}\tfrac32 x^2&\text{for }0\leq x \leq \tfrac12\\
x-\tfrac12\,(1-x)^2&\text{for }\tfrac12\leq x \leq 1,
\end{cases}$$
which indeed lies in between $F_{U_1}(x)=x$ and $F_{\max\{U_1,U_2\}}(x)=x^2$, see Figure \ref{plot}.
\begin{figure}[H]
	\centering
	\includegraphics[scale=0.75]{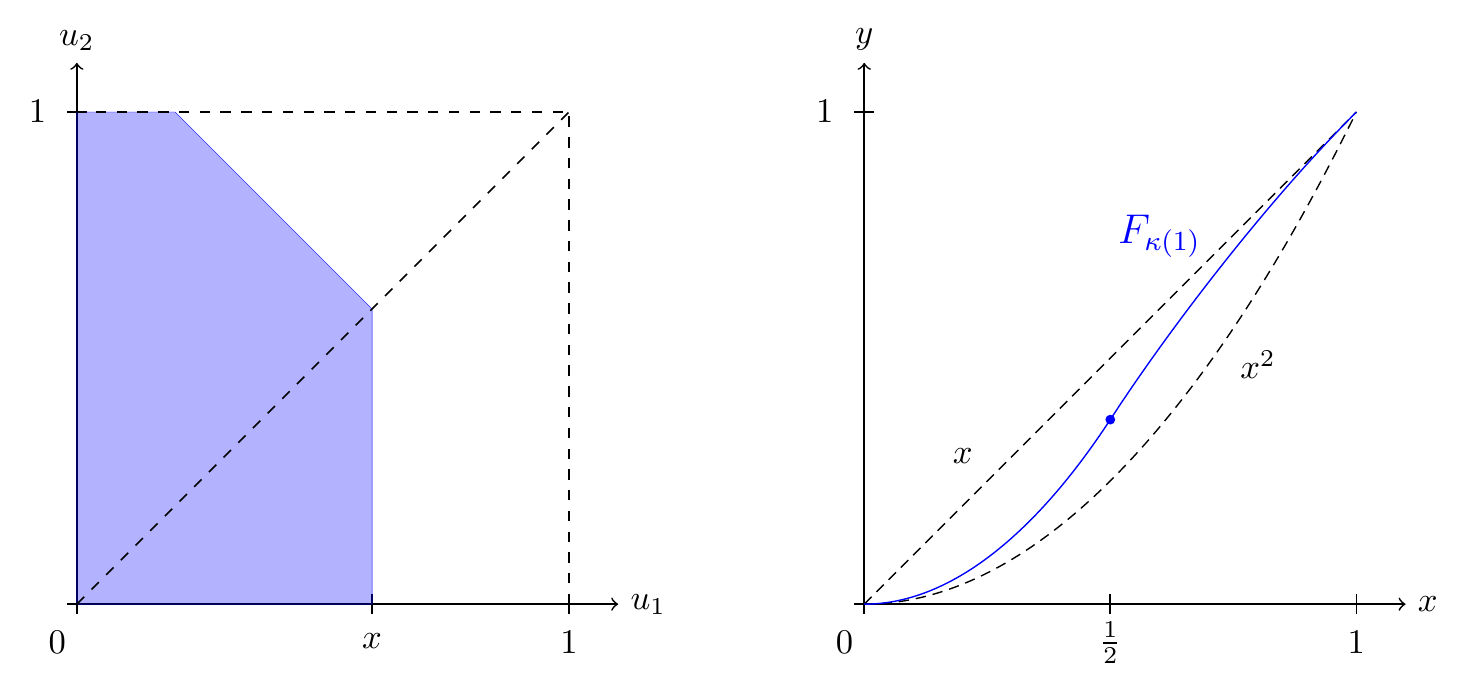}
	\caption{On the left a visualization of $\Prob(\kappa_G(1)\leq x)$, on the right the distribution
		function of $\kappa_G(1)$.\label{plot}}
\end{figure}
\item The simplest nontransitive graph is the path on three vertices, i.e.

\vspace*{-1.5em}
\makebox[\textwidth][c]{
	\begin{minipage}[t]{\textwidth}
		\begingroup
		\begin{minipage}[t]{0.45\textwidth}
			 $$G=(\{1,2,3\},\{\langle1,2\rangle, \langle 2,3\rangle\}).$$
		\end{minipage}
		\hfill
		\begin{minipage}[t]{0.45\textwidth}
			\begin{figure}[H]
				\includegraphics[scale=1]{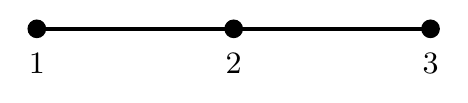}
			\end{figure}
		\end{minipage}
		\endgroup
\end{minipage}}

Simple conclusions based on Lemma \ref{collection} reveal that the suprema of achievable water levels are given by the random variables
\begin{align*}
\kappa_{G}(1)&=\max\big\{\eta_0(1),\tfrac{\eta_0(1)+\eta_0(2)}{2},
\tfrac{\eta_0(1)+\eta_0(2)+\eta_0(3)}{3}\big\},\\
\kappa_{G}(2)&=\max\big\{\eta_0(2),\tfrac{\eta_0(1)+\eta_0(2)}{2},\tfrac{\eta_0(2)+\eta_0(3)}{2},
\tfrac{\eta_0(1)}{2}+\tfrac{\eta_0(2)+\eta_0(3)}{4},\\
&\hspace{1.35cm}\tfrac{\eta_0(3)}{2}+\tfrac{\eta_0(1)+\eta_0(2)}{4}\big\}.
\end{align*}

For independent $\text{unif}(0,1)$ initial water levels, a tedious but elementary calculation
reveals the (already somewhat involved) distribution function
\begin{align*}
F_{\kappa_G(1)}(x) 
        &=\begin{cases}\tfrac83 x^3\\
                       -\tfrac{11}{6}x^3+\tfrac92 x^2-\tfrac32 x+\tfrac16\\
                       -\tfrac{23}{6}x^3+\tfrac{13}{2} x^2-2x+\tfrac16\\
                       \tfrac{2}{3}x^3-\tfrac52 x^2+4x-\tfrac76\\
          \end{cases}\text{ for }
          \begin{cases}x \in [0,\tfrac13]\\x \in [\tfrac13,\tfrac12]\\ 
                       x \in [\tfrac12,\tfrac23]\\ x \in [\tfrac23,1].\\
          \end{cases}
\end{align*}
\end{enumerate}
\end{example}

Note that in general for a finite graph $G=(V,E)$, we trivially get
$$\eta_0(v)\leq\kappa_G(v)\leq \max\{\eta_0(u),\;u\in V\},$$ consequently
$\supp\big(\kappa_G(v)\big)=[0,1]$, just the same as for the average.
As mentioned above, when dealing with infinite graphs, however, $\kappa_G(v)$ shows a different behavior
than the global average: In order to determine whether the supremum of achievable water levels at a
given vertex $v$ is a.s.\ constant or not, we have to investigate the structure of the infinite
graph a bit more closely.

What turns out to be crucial is, whether the graph contains an infinite self-avoiding path with sufficiently many extra neighbors as a subgraph. If it does, the distribution of
$\kappa_G(v)$ becomes degenerate for all $v\in V$, see Theorem \ref{inf}: One
can in fact, with probability $1$, push the water level at $v$ arbitrarily close to 1, the essential supremum of the marginal distribution.

The two-sided infinite path, however, is too lean to feature such a substructure and behaves
therefore much more like a finite graph, in the sense that the distribution of $\kappa_G(v)$ is nondegenerate -- see Theorem \ref{qtgraphs}. In order to develop these two results, let us first
properly define what we mean by ``sufficiently many extra neighbors''.

\begin{figure}[H]
	\centering
	\includegraphics[scale=0.88]{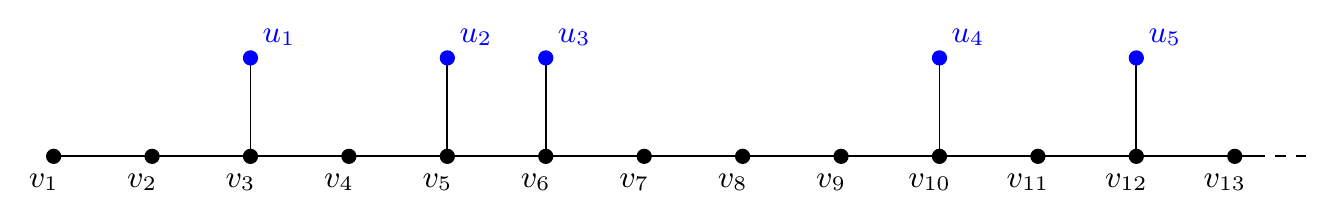}
	\caption{The beginning part of a neighbor-rich half-line.\label{half_line}}
\end{figure}

\begin{definition}\label{hl}
Let $G=(V,E)$ be an infinite, connected, simple graph. It is said to contain a {\em neighbor-rich half-line},
if there exists a subgraph of $G$ consisting of a half-line
$$H=\big(\{v_n,\;n\in\N\},\{\langle v_n, v_{n+1}\rangle,\;n\in\N\}\big)$$
and distinct vertices $\{u_n,\;n\in\N\}$ from $V\setminus\{v_n,\;n\in\N\}$ such that there is an
injective function $f:\N\to\N$ with the following two properties (cf.\ Figure \ref{half_line}):
\begin{enumerate}[(i)]
\item For all $n\in\N$: $\langle u_n, v_{f(n)}\rangle\in E$, i.e.\ the vertices $u_n$ and $v_{f(n)}$ are
neighbors in $G$.
\item The function $f$ is growing slowly, in the sense that $\sum_{n=1}^\infty\frac{1}{f(n)}$ diverges.
\end{enumerate}
\end{definition}

\noindent
Note that we can always take the function $f$ to be (strictly)
increasing -- by relabeling $\{u_n, \;n\in\N\}$ if necessary. Furthermore, if $G$ is connected and contains a neighbor-rich half-line, we can choose any vertex
$v\in V$ to be its beginning vertex: If $v_k$ is the vertex with highest index at shortest distance to $v$ in
$H$, replace $(v_1,\dots,v_k)$ by a shortest path from $v$ to $v_k$ in $H$. The altered half-line will still
be neighbor-rich, since for all $M,N\in\N$ and $f$ as above:
$$\sum_{n=1}^\infty\frac{1}{f(n)}=\infty\quad\Longleftrightarrow\quad\sum_{n=M}^\infty\frac{1}{f(n)+N}=\infty.$$

\noindent
With this notion in hand, we can state and prove the following result:

\begin{theorem}\label{inf}
Consider an infinite (connected) graph $G=(V,E)$ and the initial water levels to be i.i.d.\ $\text{\upshape
unif}(0,1)$. Let $v\in V$ be a fixed vertex of the graph.
If $G$ contains a neighbor-rich half-line, then $\kappa_G(v)=1$ almost surely.
\end{theorem}

Before we get down to the proof of this theorem, let us provide a standard auxiliary result which will come in useful:

\begin{lemma}\label{diverg}
For $\epsilon>0$, let $(Y_n)_{n\in\N}$ be an i.i.d.\ sequence having Bernoulli distribution with parameter
$\epsilon$. If the function $f:\N\to\N$ is strictly increasing and such that $\sum_{n=1}^\infty\frac{1}{f(n)}$
diverges, then
$$\sum_{n=1}^\infty \frac{Y_n}{f(n)}=\infty\quad\text{almost surely.}$$
\end{lemma}

\begin{proof}
Let us define
$$X_k=\sum_{n=1}^k \frac{Y_n-\epsilon}{f(n)}\quad \text{for all }k\in\N.$$
As the increments are independent, bounded and centered, this defines a martingale with respect to the natural filtration.
Furthermore,
$$\E(X_k^2)=\sum_{n=1}^k \frac{\E(Y_n-\epsilon)^2}{f(n)^2}=(\epsilon-\epsilon^2)\cdot\sum_{n=1}^k\frac{1}{f(n)^2}
\leq\epsilon\,\frac{\pi^2}{6}.$$
By the $L^p$-convergence theorem (see for instance Thm.\ 5.4.5 in \cite{Durrett}) there exists a random variable $X$
such that $X_k$ converges to $X$ almost surely and in $L^2$. Having finite variance, $X$ must be a.s.\ real-valued
and due to
$$\sum_{n=1}^k \frac{Y_n}{f(n)}-X_k=\epsilon\cdot\sum_{n=1}^n \frac{1}{f(n)},$$
the divergence of $\sum_{n=1}^\infty\frac{1}{f(n)}$ forces
$\sum_{n=1}^\infty \frac{Y_n}{f(n)}=\infty$ almost surely.
\end{proof}\vspace*{1em}

\begin{nproof}{of Theorem \ref{inf}}	
Given a graph $G$ with the properties stated and a vertex $v$, we can choose a neighbor-rich half-line $H$ with
$v=v_1$ and the set of extra neighbors $\{u_n{}\}_{n\in\N}$ as described in and after Definition \ref{hl}.
The initial water levels at $\{u_n\}_{n\in\N}$ are i.i.d.\ $\mathrm{unif}(0,1)$, of course.

Depending on the random initial profile, let us define the following SAD-procedure starting at $v$:
Fix $\epsilon,\delta>0$ and let $\{N_k\}_{k\in\N}$ be the increasing (random) sequence of indices, chosen such that the initial water level at $u_{N_k}$ is at least $1-\epsilon$ for all $l$. Then define the SAD-procedure --
starting with $\xi_0(v)=1,\ \xi_0(u)=0$ for all $u\in V\setminus\{v\}$ -- such that
first all vertices along the path $(v_1,v_2,\dots,v_{f(N_1)},u_{N_1})$ exchange liquids sufficiently often to get
$$\xi_{t_1}(u_{N_1})\geq\frac{1}{f(N_1)+2}\quad\text{for some }t_1>0,$$
and never touch $u_{N_1}$ again. Note that by Lemma \ref{evenout}, $\xi_{k}(u_{N_1})$ can be pushed
as close to $\tfrac{1}{f(N_1)+1}$ as desired in this way.
At time $t_1$, the joint amount of water in the glasses at $v_1,v_2,\dots,v_{f(N_1)}$ equals $1-\xi_{t_1}(u_{N_1})$ and
we will repeat the same procedure along $(v_1,v_2,\dots,v_{f(N_2)},u_{N_2})$ to get
$$\xi_{t_2}(u_{N_2})\geq\frac{1}{f(N_2)+2}\cdot\big(1-\xi_{t_1}(u_{N_1})\big)\quad\text{for some }t_2>t_1$$
and iterate this.

After $m$ iterations of this kind, the joint amount of water localized at vertices of the half-line $H$ equals
$1-\sum_{k=1}^m \xi_{t_k}(u_{N_k})$, which using $1-x\leq\text{e}^{-x}$ can be bounded from above as follows:
\begin{align}\label{uppbound}\begin{split}
1-\sum_{k=1}^m \xi_{t_k}(u_{N_k})&\leq\prod_{k=1}^m \bigg(1-\frac{1}{f(N_k)+2}\bigg)\\
                               &\leq \exp\bigg(-\sum_{k=1}^m\frac{1}{f(N_k)+2}\bigg).\end{split}
\end{align}
Defining $Y_n:=\mathbbm{1}_{\{\eta_0(u_n)\geq 1-\epsilon\}}$ for all $n\in\N$ we get $(Y_n)_{n\in\N}$ i.i.d.\ 
$\text{Ber}(\epsilon)$ and can rewrite the limit of the sum in the exponent as follows:
$$\sum_{k=1}^\infty\frac{1}{f(N_k)+2}=\sum_{n=1}^\infty\frac{Y_n}{f(n)+2}.$$
This allows us to conclude from Lemma \ref{diverg} that the exponent in (\ref{uppbound}) tends a.s.\ to $-\infty$
as $m \to \infty$. Consequently, $m,T\in \N$ can be chosen large enough such that with probability $1-\delta$ it holds that
$$\sum_{k=1}^m \xi_{t_k}(u_{N_k})\geq 1-\epsilon\quad\text{and}\quad t_m\leq T.$$

Given this event, the move sequence corresponding to the SAD-procedure just described -- adding no further
updates after time $t_m$, i.e.\ $\mu_t=0$ for $t_m<t\leq T$ -- then ensures
(cf.\ Lemma \ref{dual}) that
$$\eta_T(v)\geq\sum_{k=1}^m\xi_{T}(u_{N_k})\,\eta_0(u_{N_k})\geq(1-\epsilon)^2,$$
forcing $\kappa_G(v)\geq(1-\epsilon)^2$ with probability at least $1-\delta$. Since $\delta>0$ was arbitrary, this implies
$\kappa_G(v)\geq(1-\epsilon)^2$ a.s.\ and letting $\epsilon$ go to $0$ then establishes the claim.
\end{nproof}

\begin{nproof}{of Theorem \ref{qtgraphs}}	
	Given i.i.d.\ $\mathrm{unif}(0,1)$ initial water levels on an infinite (connected) graph $G=(V,E)$, the strong law of large numbers guarantees $\kappa_G(v)\geq\tfrac12$ (cf. Lemma \ref{evenout}). In
	fact, we can even get the strict inequality: Take $\{v_k,\;k\in\N\}\subseteq V$ such that $v_1=v$
	and for all	$n\in\N$, $\{v_1,\dots,v_n\}$ is a connected vertex-set of size $n$.
	Further, define $X_k=\eta_0(v_k)-\frac12$ for all $k\in\N$.
	Then $(S_n)_{n\in\N}$, defined by $S_n=\sum_{k=1}^n X_k$, is an i.i.d.\ symmetric random walk, whence
	$\limsup_{n\to\infty} S_n=\infty$ almost surely, see e.g.\ Thm.\ 4.1.2 in \cite{Durrett}. In particular, there	exists a.s.\ an $n\in\N$ with $\frac1n\cdot\sum_{k=1}^n \eta_0(v_k)>\frac12$. Again by Lemma \ref{evenout}, this ensures $\Prob\big(\kappa_G(v)>\tfrac12\big)=1$.
	
	If $G$ is the two-sided infinite path, there is a positive probability that the vertex $v$ is what Häggström \cite{ShareDrink}
	calls two-sidedly $\epsilon$-flat with respect to the initial profile (cf.\ Lemma 4.3 in \cite{ShareDrink}), i.e.
	\begin{equation}\label{tflat}
	\frac{1}{m+n+1}\sum_{u=v-m}^{v+n}\eta_0(u)\in\left[\tfrac12-\epsilon,\tfrac12+\epsilon\right]\quad
	\text{for all }m,n\in\N_0.
	\end{equation}	
	Lemma 6.3 in \cite{ShareDrink} states that in this situation, the water level at $v$ is bound to
	stay within the interval $[\tfrac12-6\epsilon,\tfrac12+6\epsilon]$, irrespectively of future updates.
	Together with the simple observation that $\kappa_G(v)\geq\eta_0(v)$ in general, this implies $\{\frac12,1\}\subseteq\mathrm{supp}\big(\kappa_\Z(v)\big)$.
	
	In order to establish the first part, we are left to show that $\kappa_\Z(v)<1$ almost surely.
	Let us assume $\Prob\big(\kappa_\Z(v)=1\big)>0$ for contradiction and fix a finite set $I\subseteq\Z$ of nodes. If $\kappa_\Z(v)=1$, there exist finite move sequences achieving water levels at $v$ arbitrarily close to 1. Since $\max\{\eta_0(u),\; u\in I\}<1$ a.s., the corresponding dual SAD-profiles must tend to 0 on $I$ (cf.\ Lemma \ref{dual}).
	As a consequence, $\mathbbm{1}_{\{\kappa_\Z(v)=1\}}$ is independent of $\{\eta_0(u),\;u\in I\}$
	and therefore $\{\kappa_\Z(v)=1\}$ a tail event. By Kolmogorov's zero-one law and our assumption,
	it has to be an almost sure event, which contradicts
	$\frac12\in\mathrm{supp}\big(\kappa_\Z(v)\big)$.	
	
	In view of Theorem \ref{inf}, to prove the second part, we only have to verify, that an infinite, connected,
	quasi-transitive graph that is not the two-sided infinite path contains a neighbor-rich half-line. Since $G$ is infinite (and by our
	general assumptions both connected and having finite maximal degree) a compactness argument guarantees the existence of
	a half-line $H$ on the vertices $\{v_n,\;n\in\N\}$ such that $v_1=v$ and the graph distance from $v_n$ to $v$ is $n-1$ for
    all $n$.
	
	Let us consider the function $d:V\to\N_0$, where $d(u)$ is the graph distance from the node $u$ to a vertex of degree
	at least 3 being closest to it.
	Since $G$ is quasi-transitive, connected and not the two-sided infinite path, $d$ is finite and can take on only
	finitely many values, which is why it has to be bounded, by a constant $c\in \N$ say. Consequently, $G$ can not
	contain stretches of
	more than $2c$ linked vertices of degree 2. For this reason, there must be a vertex among $v_3,\dots,v_{2c+3}$,
	say $v_{f(1)}$, having a neighbor $u_1$ outside of $H$. In the same way, we can find a vertex $u_2$ outside $H$
	having a neighbor $v_{f(2)}$ among $v_{2c+6},\dots,v_{4c+6}$ and in general some $u_n$ not part of
	$H$ but linked to a vertex $v_{f(n)}\in\{v_n,\;n\in\N\}$ with $$(n-1)\,(2c+3)+3\leq f(n)\leq n\,(2c+3)\quad\text{for all }n\in\N.$$
	This choice ensures that $v_{f(k)}$ and $v_{f(n)}$ are at graph distance at least 3 for $k\neq n$, which forces the
	set $\{u_n,\;n\in\N\}$ to consist of distinct vertices.	Due to
	$$\sum_{n=1}^\infty\frac{1}{f(n)}\geq\frac{1}{2c+3}\,\sum_{n=1}^\infty\frac{1}{n}=\infty,$$
	$H$ is a neighbor-rich half-line in the sense of Definition \ref{hl} as desired, which concludes the proof.
\end{nproof}

\begin{remark}
\begin{enumerate}[(a)]
\item
Note that the essential property of the initial water levels, needed in the proof of Theorem \ref{inf}, was
independence. The argument can immediately be generalized to the situation where the initial water levels are independently (but not
necessarily identically) distributed on a bounded interval $[0,C]$ and we have some weak form of uniformity, namely:

For every $\delta>0$, there exists some $\epsilon>0$ such that
for all $v\in V$: $$\Prob\big(\eta_0(v)>C-\delta\big)\geq \epsilon.$$
The sequence $Y_n:=\mathbbm{1}_{\{\eta_0(u_n)\geq C-\delta\}},\ n\in\N,$ similar to the one defined in the proof of
Theorem \ref{inf} will no longer be i.i.d.\ $\text{Ber}(\epsilon)$, but an appropriate coupling will ensure that
$$\sum_{n=1}^\infty \frac{Y_n}{n}\geq\sum_{n=1}^\infty \frac{Z_n}{n}\quad\text{a.s.,}$$
where $(Z_n)_{n\in\N}$ is an i.i.d.\ sequence of $\text{Ber}(\epsilon)$ random variables. Accordingly,
we get $\kappa_G(v)=C$ a.s.\ even in this generalized setting.

\item To see how surprising Corollary \ref{main} is, consider the case of the quasi-transitive graph
$G_M$, obtained by taking $\Z$ and adding a single twig (an edge leading to a single vertex and nothing
more) to every $M$th vertex of $\Z$. Intuitively, one would think that for large $M$, e.g.\ 
$M=10^{100}$, the distribution of $\kappa_{G_M}(v)$ at a vertex $v$ far from any of the twigs would be
qualitatively similar to the nondegenerate distribution of $\kappa_\Z(v)$. But this is not the case,
since according to Theorem \ref{qtgraphs}, the distribution of $\kappa_{G_M}(v)$ is concentrated at 1.

\item
Note that the second part of Theorem \ref{qtgraphs} trivially implies $\kappa_G(v)=1$ a.s.\ 
for any graph $G$ containing an infinite, connected, quasi-transitive subgraph.
It is worth emphasizing that Theorem \ref{qtgraphs} does, however, not capture the full statement of Theorem \ref{inf}:
If we take $G$ to consist of the two-sided infinite path $\Z$ and a single twig added to every node that corresponds to a prime number,
the only infinite, connected, quasi-transitive subgraph of $G$ is $\Z$ itself. However, since $G$ contains a neighbor-rich
half-line, Theorem \ref{inf} states that $\kappa_G(v)=1$ a.s.\ for i.i.d.\ $\mathrm{unif}(0,1)$ initial water levels
and any target vertex $v$.

\item\label{concl}
As alluded to in the introduction, the results from Section \ref{infinite} do have implications for
the Deffuant model for opinion formation, whose analysis originally inspired the optimization problem
investigated here. The qualitative change of behavior of $\kappa_G(v)$ for a fixed vertex $v\in\Z^d$, when
switching from the two-sided infinite path ($d=1$) to higher-dimensional grids ($d\geq2$), constitutes
a substantial hinder in the analysis of the Deffuant model: As proved in Lemma 6.3 in \cite{ShareDrink},
for i.i.d.\ $\mathrm{unif}(0,1)$ initial values on $\Z$, almost surely there exist nodes that are stuck
with a value close to the mean of the marginal distribution (irrespectively of future interactions of
neighboring nodes), while Theorem \ref{qtgraphs} tells us that this is no longer true in higher
dimensions. These vertices played however a central role in the two analyses of the Deffuant model on
$\Z$, done by Lanchier \cite{Lanchier} and Häggström \cite{ShareDrink} respectively, which is why their
ideas can not be transferred to the case $d\geq2$ immediately.
\end{enumerate}
\end{remark}



\vspace{0.5cm}
\makebox[\textwidth][c]{
\begin{minipage}[t]{1.2\textwidth}
\begingroup
	\begin{minipage}[t]{0.5\textwidth}
	{\sc \small Olle Häggström\\
   Department of Mathematical Sciences,\\
   Chalmers University of Technology,\\
   412 96 Gothenburg, Sweden.}\\
   olleh@chalmers.se
	\end{minipage}
	\hfill
	\begin{minipage}[t]{0.5\textwidth}
	{\sc \small Timo Hirscher\\
   Institute for Mathematics,\\
   J.W.\ Goethe University,\\
   60054 Frankfurt a.M., Germany.}\\
   hirscher@math.uni-frankfurt.de\\
	\end{minipage}
	\endgroup
\end{minipage}}


\begin{thebibliography}{9}
    \bibitem{Durrett}
         {\sc Durrett, R.,}
         {\em ``Probability: Theory and Examples (4th edition)''},
         Cambridge University Press, 2010.
    \bibitem{ShareDrink}
         {\sc Häggström, O.,}
         {A pairwise averaging procedure with application to consensus formation in the Deffuant model,}
         {\em Acta Applicandae Mathematicae},
         Vol.\ 119 (1), pp.\ 185-201, 2012.
    \bibitem{Deffuant}
         {\sc Häggström, O.} and {\sc Hirscher, T.,}
         {Further results on consensus formation in the Deffuant model,}
         {\em Electronic Journal of Probability},
         Vol.\ 19, 2014.
    \bibitem{Lanchier}
         {\sc Lanchier, N.,}
         {The critical value of the Deffuant model equals one half,}
         {\em Latin American Journal of Probability and Mathematical Statistics}, 
         Vol.\ 9 (2), pp.\ 383-402, 2012.
    \bibitem{Shang}
         {\sc Shang, Y.,}
         {Deffuant model with general opinion distributions: First impression and critical confidence bound,}
         {\em Complexity},
         Vol.\ 19 (2), pp.\ 38-49, 2013.
 \end{thebibliography}
\end{document}